\documentclass[12pt]{article}
\usepackage{latexsym, amsmath, amssymb}
\def\dim{\mathop{\hbox{dim}}}

\def\PI{\mathop{\hbox{Im}}}

\def\Ric{\mathop{\hbox{Ric}}}
\def\Ke{\mathop{\hbox{K}}}

\newtheorem{theorem}{Theorem}[section]
\newtheorem{proposition}[theorem]{Proposition}
\newtheorem{remark}[theorem]{Remark}

\newtheorem{example}[theorem]{Example}

\def\Im{\mbox{Im}}
\def\Re{\mbox{Re}}

\def\ds{\displaystyle}

\def\veps{\varepsilon}
\def\1{\mbox{\bf 1}}

\def\R{\mathbb R}
\def\Z{\mathbb Z}
\def\N{\mathbb N}
\def\C{\mathbb C}

\def\1{\mathbb 1}

\def\al{\alpha}

\newtheorem{teor}{Theorem}[section]

\newtheorem{remar}[teor]{Remark}

\newtheorem{corol}[teor]{Corollary}
\newtheorem{lemma}[teor]{Lemma}

\newcommand{\fdim}{\hspace*{\fill}$\Box$}
\newcommand{\dimostr}{{\bf Proof: }}

\newcommand{\real}{\Bbb{R}}

\newcommand{\complex}{\Bbb{C}}

\begin{document}
\centerline{\large\bf TYZ expansion for the Kepler manifold 
\footnote{
The first author was supported in part by GNAMPA-INDAM, Italy.
The second author was supported in part by the M.I.U.R.
Project \lq\lq Geometric Properties of Real and Complex
Manifolds''.}}



\vspace{0.5cm}

\centerline{\small Todor Gramchev and Andrea Loi}
\centerline{\small Dipartimento di Matematica e Informatica
-- Universit\`{a} di Cagliari}
\centerline{\small Via Ospedale 72, 09124 Cagliari-- Italy}
\centerline{\small e-mail : todor@unica.it, loi@unica.it}

\vskip 0.5cm

\centerline{\bf Abstract}
{\small

The main goal of the paper is to address the issue of  the existence 
of Kempf's distortion function and the Tian-Yau-Zelditch (TYZ) asymptotic expansion for the Kepler manifold - an important example 
of non compact manfold. 
Motivated by the recent results for compact manifolds
we construct Kempf's distortion function  and derive a precise  TYZ asymptotic expansion 
for the Kepler manifold. 
We get an exact formula:  finite asymptotic expansion of $n-1$ terms and  
exponentially small error terms  uniformly with respect to the discrete quantization parameter $m\rightarrow \infty $ 
and  $\rho \rightarrow \infty$, $\rho$ being the polar radius in $\C^n$. 
 Moreover,  
the coefficents are calculated explicitly and they turned out to be  homogeneous functions with respect to the polar radius in the Kepler manifold. 
We also prove and derive an asymptotic expansion of  the  obtstruction  term with the coefficients being defined by geometrical quantities. 
We  show that our estimates are sharp by  analyzing the nonharmonic behaviour of $T_m$ and the error term of the approximation of the Fubini--Study metric by $m\omega$ for $m\rightarrow +\infty$. 
The arguments of the proofs combine geometrical methods, quantization tools  and  
functional analytic techniques for investigating  asymptotic expansions in the framework of analytic-Gevrey spaces. 
\vskip 0.3cm

\noindent
{\it{Keywords}}: K\"{a}hler manifolds; quantization; quantum mechanics;
TYZ asymptotic expansion; exponential reminder.

\noindent
{\it{Subj.Class}}: 53C55, 58F06, 58J37 }

\section{Introduction and statements of the main results}
Let $g$ be  a K\"{a}hler metric
on a  complex $n$-dimensional manifold  $M$.
Assume that $g$ is polarized with respect to
a holomorphic line bundle  $L$ over $M$,
i.e. $c_1(L)=[ \omega ]$, where $\omega$
is the K\"{a}hler form associated to $g$
and $c_1(L)$ denote the first Chern class of $L$.
Let $m\geq 1$ be a non-negative integer and
let $h_m$ be an  Hermitian metric
on $L^m=L^{\otimes m}$ such that its Ricci curvature $\Ric (h_m)=m\omega$.
Here $\Ric (h_m)$ is the two form on $M$
whose local expression is given by
\begin{equation}\label{rich}
\Ric (h_m)=-\frac{i}{2}
\partial\bar\partial\log h_m(\sigma (x), \sigma (x)),
\end{equation}
for a trivializing holomorphic
section $\sigma :U\rightarrow L^m\setminus\{0\}$.
In the quantum mechanics terminology
$L^m$ is called the {\em quantum line bundle},
 the pair $(L^m, h_m)$ is called a {\em geometric
quantization}
of the K\"{a}hler manifold $(M, m\omega)$
and 
$h = m^{-1}$ play the role of Planck's constant
(see e.g. \cite{arlquant}).
Consider the separable complex
Hilbert space ${\cal H}_m$
consisting of global holomorphic
sections  $s$ of $L^m$ such that
$$\langle s, s\rangle_m=
\int_Mh_m(s(x), s(x))\frac{\omega^n}{n!}<\infty .$$

Let $x\in M$ and
$q\in L^m \setminus \{0\}$ a fixed point of the fiber over $x$.
If one evaluates
$s\in {\cal H}_m$ at $x$,
one gets a multiple $\delta_{q}(s)$ of $q$,
i.e. $s(x)=\delta_{q}(s)q$.
The map $\delta_{q}:{\cal H}_m\rightarrow {\complex}$
is a continuous linear functional
\cite{cgra}. Hence from Riesz's
theorem,
there exists a unique $e_{q}^m\in {\cal H}$
such that
$\delta_{q}(s)=\langle s, e_{q}^m\rangle_m, \forall s\in {\cal H}_m$, i.e.
\begin{equation}\label{coherentstates}
s(x)=\langle s, e_{q}^m\rangle_mq .
\end{equation}
It follows that
$$e_{cq}^m=
\overline{c}^{-1}e_{q}^m,\    \forall c\in {\complex}^*.$$
The holomorphic section $e_{q}^m\in {\cal H}_m$ is
called the {\em coherent state} relative to the point $q$.
Thus, one can define
a smooth function on $M$
\begin{equation}\label{defep}
T_{m} (x) =h_m(q, q)\|e_{q}^m\|^2,\ \|e_{q}^m\|^2=\langle e_q^m,
e_{q}^m\rangle,
\end{equation}
where $q\in L^m\setminus \{0\}$ is any point
on the fiber of $x$.
If $s_j,\,\, j=0,\dots d_m,$ ($d_m+1=\dim {\cal H}_m\leq\infty$) is a
orthonormal basis for
$({\cal H}_m, \langle\cdot ,\cdot\rangle_m)$ then one can easily verify that
\begin{equation}\label{Tmo}
T_{m} (x) =\sum_{j=0}^{d_m}h_m(s_j(x), s_j(x)).
\end{equation}

Notice that when $M$ is compact
${\cal H}_m=H^0(L^m)$, where $H^0(L^m)$
denotes the space of global holomorphic sections
of $L^m$. Hence in this case $d_m<\infty$
and (\ref{Tmo}) is a finite sum.

The function $T_{m}$ has appeared in the
literature under different names.
The earliest one was probably the $\eta$-function
of J. Rawnsley \cite{ra1} (later renamed to $\epsilon$
function in \cite{cgra}), defined for arbitrary K\"{a}hler manifolds,
followed by the {\em distortion function }
of Kempf \cite{ke} and Ji \cite{ji},
for the special case of Abelian varieties
and of Zhang \cite{zha} for complex projective varieties.
The metrics for which $T_m$
is constant were called {\em critical} in \cite{zha}
and {\em balanced} in \cite{do}
(see also \cite{arlcomm}, \cite{bergbal},
\cite{baC} and  \cite{regcov}).
If $T_m$ are constants for all sufficiently
large $m$ then the  geometric quantization
$(L^m, h_m)$
associated to the K\"{a}hler manifold $(M, g)$
is called regular.
Regular quantization play a prominent role in the theory of
quantization by deformation of K\"{a}hler manifolds developed
in \cite{cgra} (see also \cite{regscal}). 

Fix $m\geq 1$.
Under the hypothesis that for
each point $x\in M$
there exists $s\in {\cal H}_m$
non-vanishing at $x$,
one can give
a geometric interpretation
of $T_m$   as follows.
Consider the holomorphic map
of $M$ into the complex
projective space
${\complex}P^{d_m}$:
\begin{equation}\label{psiglob}
\varphi_m :M\rightarrow {\complex}P^{d_m}:
x\mapsto [s_0(x): \dots :s_{d_m}(x)].
\end{equation}

One can prove that
\begin{equation}\label{obstr}
\varphi ^*_m(\omega_{FS})=m\omega +
\frac{i}{2}\partial\bar\partial\log T_{m} ,
\end{equation}
where $\omega_{FS}$ is the Fubini--Study form on
${\complex}P^{d_m}$, namely the form which
in homogeneous coordinates
$[Z_0,\dots, Z_{d_m}]$ reads as
$\omega_{FS}=\frac{i}{2}\partial\bar\partial\log \sum_{j=0}^{d_m}
|Z_j|^2$.

Clearly \eqref{obstr} leads to 
\begin{equation}\label{obstrer1} 
\frac{\varphi_m ^*(\omega_{FS})}{m}- \omega =
\frac{i}{2m}\partial\bar\partial\log T_{m} ,
\end{equation}
therefore the term 
\begin{equation}\label{obstrer2} 
{\mathcal E}_m(x) :=
\frac{i}{2m}\partial\bar\partial\log T_{m} ,
\end{equation}
turns out to play a role of the ``error'' of the approximation of $ \omega$ (resp. $g$) by 
$\frac{\varphi_m ^*(\omega_{FS})}{m}$ (resp. $\frac{\varphi_m ^*(g_{FS})}{m}$). 

\vspace{10pt}

Observe that by (\ref{obstr}),
if there exists  $m$ such that  $mg$ is a balanced metric, or more generally
if $T_{m}$ is harmonic, then ${\mathcal E}_m(x)$ is identically zero and
hence $mg$ is projectively induced
via the coherent states map
$\varphi_m$ (see \cite{arlquant} for more details
on the link between 
projectively induced K\"{a}hler metrics
and balanced metrics).
Recall that a  K\"{a}hler metric $g$ on a complex manifold
$M$ is  {\em projectively induced} if there exists a K\"{a}hler  (i.e. a holomorphic and isometric)
immersion $\psi : M\rightarrow {\complex}P^N, N\leq \infty$
such that $\psi^*(g_{FS})=g$.
Projectively induced K\"{a}hler metrics enjoyes very nice properties 
and they were deeply studied in \cite{ca} (see also the begining of Section \ref{sharp} below).
Not all K\"{a}hler metrics
are balanced or  projectively induced.
Nevertheless,
when $M$ is compact, Tian \cite{ti0}
and Ruan \cite{ru}
solved a conjecture posed by Yau by proving
that the sequence of metrics
$\frac{\varphi_m ^*(\omega_{FS})}{m}$
$C^{\infty}$-converges to $\omega$.
In other words, any polarized metric on 
compact complex manifold is,
the $C^{\infty}$-limit of
(normalized) projectively induced 
K\"{a}hler metrics.
Zelditch \cite{ze}
generalized Tian--Ruan theorem
by proving 
a complete asymptotic expansion in the $C^\infty$ category, namely  
\begin{equation}\label{asymptoticZ}
T_{m}(x) \sim 
\sum_{j=0}^\infty  a_j(x)m^{n-j}
\end{equation} 
where  $a_j$, $j=0,1, \ldots$, are smooth coefficients with $a_0(x)=1$, and 
for any nonnegative integers $r,k$ the following estimates hold: 
\begin{equation}\label{rest}
||T_{m}(x)-
\sum_{j=0}^{k}a_j(x)m^{n-j}||_{C^r}\leq C_{k, r}m^{n-k-1},
\end{equation}
where $C_{k, r}$ are constant depending on $k, r$ and on the K\"{a}hler
form $\omega$ and $ || \cdot ||_{C^r}$ denotes  the $C^r$ norm in local coordinates..

Later on, Lu \cite{lu}, by means of  Tian's peak section method,
proved  that 
each of the coefficients $a_j(x)$ in 
(\ref{asymptoticZ}) is a polynomial of the curvature and its covariant
derivatives at $x$ of the metric $g$. 
Such a polynomials can be found by finitely many steps
of algebraic operations. Furthermore $a_1(x)=\frac{1}{2}\rho$,
where $\rho$ is the scalar curvature of the polarized metric $g$
(see also \cite{loianal} and \cite{loismooth} for
the computations of the coefficients $a_j$'s through
Calabi's diastasis function). 

The expansion (\ref{asymptoticZ}) is called
the  {\em TYZ (Tian--Yau--Zelditch) expansion}.
\vspace{10pt}

The aim of the present paper is  to adress the problem of  
TYZ expansions for noncompact manifolds. Our motivations is twofolded. First,
its is purely geometrical question of its own interest. 
Secondly, 
 we are inspired by the previous 
works  of M. Engli\v{s} \cite{engber}, \cite{enpseudo}, 
\cite{me2}, where analytical tools from the the theory of asymptotic expansions  
have been applied in order to extend
 Berezin's quantization method cf. \cite{Ber1},
 \cite{Ber2}
 to non homogeneous complex domains on ${\complex}^n$ (see also 
 \cite{Odz1}, \cite{Odz2}, \cite{MlaTsa}).




We choose as a noncompact manifold  the Kepler manifold $(X, \omega)$,
namely the cotangent bundle of the $n$-dimensional sphere
minus its zero section endowed with the standard symplectic 
 form $\omega$ (see   \cite{sou} and \cite{ra}). 
  This manifold has been considered by different
 authors and we bilieve our results  can be of some interest
 both from the mathematical and physical
 point of view.

We summarize the main novelties of our work: 
First, we compute explicitly 
the Kempf distortion function $T_m(x)$ 
for the Kepler manifold $(X, \omega)$.  
Secondly, based on this computation 
we find  an analogue  of Zelditch
and Lu's theorems above for $(X, \omega)$. 
More precisely, building upon the explicit representation of $T_m$ as an action of "singular derivatives" 
and using precise analytical methods pertinent to the study of nonlinear compositions in functional spaces, we  show that the TYZ expansion for the Kepler manifold has  
two remarkable features in comparison with the known results for compact manifolds: 
\begin{itemize}
\item first, the TYZ expansion is   {\em finite}. More precisely, it  consists of $n-1$ terms 
$$T_m(x) =   m^n  + \frac{(n-2)(n-1)}{2 |x|} m^{n-1} +
 \sum_{k=2}^{n-2}  \frac{2a_k}{|x|^k} m^{n-k}  + R_m(|x|),$$
 where $a_k,\  k\geq 2$ can be computed explicitely by recursive formulas.
\item secondly,  the {\em  reminder term has an
exponential small decay} $O(e^{-cm})$ as $m\rightarrow \infty$ uniformly with respect to $|x| \geq \delta>0$. 
 We point out that our exact formula   modulo exponentially small error 
for Kempf's distortion function might be viewed as an analogue to a geometric 
interpretation of exact asymptotic formulas appearing for the moment map and equivariant cohomology cf. M. F. Attiyah and R. Bott \cite{AtBot1} (see also \cite{DuHe1}, \cite{DuHe2}).
\end{itemize}
 
 We also derive  uniform analytic--Gevrey estimates for $T_m$ keeping the exponential decay 
for $m\rightarrow \infty$, $|x| \rightarrow \infty$ 
which  resemble the estimates 
in the framework of Gelfand--Shilov spaces $S^1_1$ appearing in the regularity theory 
for pseudodifferential operators cf. \cite{cgrjfa}, \cite{cgrbir}. 
Observe that as for the compact case our expansion shows that 
$g$ (the metric $g$ associated to the Kepler manifold $(M, \omega)$) 
is the $C^{\infty}$-limit of (suitable normalized) projectively induced K\"{a}hler metrics,
namely
$\lim_{m\rightarrow\infty}\frac{1}{m}\varphi_m^*(g_{FS})=g$
where
$\varphi_m: X\rightarrow {\complex}P^{\infty}$
is the coherent states map.
A geometric construction is proposed showing that our estimates are sharp. Indeed we
show that $g$ is not projectively induced, i.e.  it cannot
exist {\em any}  K\"{a}hler immersion of $(X, \omega)$
into a finite or infinite dimensional complex projective space.
The arguments use Calabi's tools which provide necessary and sufficient conditions for a K\"{a}hler metric
to be projecticely induced.
 
Finally, we investigate the asymptotic behaviour of the obstruction term  
$${\mathcal E}_m (z) = \sum_{j,\ell=1}^{n+1} {\mathcal E}_m^{j,\ell} (z) dz_j \wedge d\bar{z}_\ell $$ in \eqref{obstrer2} and prove  that the coefficients decay polynomially of the type $m^{-2}$. 
More precisely, for some $C>0$, they behave like   
\begin{equation}\label{obstrer2a} 
 \frac{C}{m^2 |z|^{3}} (1 + o(1))
\qquad m\rightarrow \infty ,
\end{equation} 
uniformly for $|z|$ away from the origin in $\C^n$.
In fact, we show an abstract theorem for the asymptotic behaviour of obstruction terms 
similar to \eqref{obstrer2} on  
conic manifolds of Kepler type. The proof is based  
 on a suitable choice of global singular coordinates parametrizing the Kepler manifold 
and the  use of implicit function theorem arguments. 
Consequently, by  (\ref{obstr}),  the metric  $g$ associated to $\omega$
can be approximated by suitable normalized 
projectively induced K\"{a}hler metrics with an error of the type $m^{-2}$, $m\rightarrow \infty$. 
\vspace{10pt}

The paper is organized as follows. We propose an explicit construction of 
the Kempf distortion function $T_m$ 
for the Kepler manifold $(X, \omega)$ in  Section \ref{kempf}. 
In Section \ref{expansion}
 we derive 
 an exact TYZ asymptotic expansion and derive the exponentially small decay for the remainder when $m\rightarrow \infty$.
In Section \ref{sharp} we prove (see Theorem \ref{mainteor2})
that our estimate is sharp.
Finally, Section \ref{sharplog}  contains the construction of 
the global singular parametrization of the Kepler manifold and 
 the study of the asymptotic behaviour of 
the logarithmic obstruction term (\ref{obstrer2}).








\section{Kempf's distortion function  for  the Kepler manifold}\label{kempf}

The (regularized) Kepler manifold \cite{sou} is
(may be identifed with)
the $2n$-dimensional
symplectic manifold
$(X, \omega )$, where
$X=T^*S^n\setminus{0}$
the cotangent bundle to the $n$-dimensional
sphere minus its zero section
endowed with the standard symplectic form
$\omega$.
This may further be identified with 
$$X=\{(e, x)\in {\real}^{n+1}\times {\real}^{n+1}|\
e\cdot e=1, x\cdot e=0, x\neq 0
\},$$
where the dot  denotes the standard scalar product
on ${\real}^{n+1}$.
In \cite{sou} J. Souriau showed that the Kepler  manifold
admits a natural complex structure. Indeed
he proved that by introducing

$$
z=|x|e+ix\in {\complex}^{n+1}= |x| (e + is), \qquad s = \frac{x}{|x|} \in S^n, $$
then $X$
is diffeomorphic to the isotropic cone
$$C=\{z\in {\complex}^{n+1}|\ z\cdot z=z_1^2+\cdots +z_{n+1}^2=0, z\neq 0
\}\subset {\complex}^{n+1}$$
and hence $X$ inherits the complex structure of
$C$ via this diffeomorphism.
Seven years later J. Rawnsley \cite{ra1} observed that the symplectic
form $\omega$ is indeed a K\"{a}hler form with respect to this complex structure
and it can be written (up to a factor) as
\begin{equation}\label{Kform}
\omega=\frac{i}{2}\partial\bar\partial |x|.
\end{equation}
Moreover, since $\omega$ is exact, it is trivially integral
and hence there exists
a holomorphic line bundle $L$ over $X$
such that $c_1(L)=[\omega]$.

For $n\geq 3$, $X$ is simply-connected so $L^m$
is  holomorphically trivial  ($L^m=X\times {\complex}$)
and  we can identify $H^0(L^m)$
with the set of  holomorphic functions of $X$.
Furthermore, we can define an Hermitian metric $h_m$
on $L^m=X\times {\complex}$ by
\begin{equation}\label{hmK}
h_m(\sigma (z), \sigma (z))=e^{-m|x|},
\end{equation}
where $\sigma :X\rightarrow X \times {\complex},$ is the global holomorphic
section such that $\sigma (z)= (z, 1)$.
It follows by (\ref{rich}) above that
the pair $(L^m, h_m)$ is indeed
a geometric quantization of
the Kepler manifold $(X, \omega)$.
Then the Hilbert space ${\cal H}_m$ consists
of the set of homorphic functions $f$ of $X$
such that
$$ \|f\|_m^2:=\int_X |f(z)|^2  e^{-m|x|} d\mu (z)<\infty,$$
where  
$$d\mu (z)=\frac{\omega ^n(z)}{n!}=(\frac{i}{2}\partial\bar\partial
|x|)^n.$$

Notice that in this case
$$T_m(z)=e^{-m|x|}K^{(m)}(z, z),$$
where $K^{(m)}(z, z)$ is the reproducing Kernel
for the Hilbert space ${\cal H}_m$.
At  p. 412 in  \cite{ra1} Rawnsley explicitly computed
$K(z, z)=K^{(1)}(z, z)$ (the reproducing kernel for ${\cal H}={\cal H}_1$)
and  hence the  corresponding   Kempf's distortion function,
which in our notations
reads as:
\begin{equation}\label{raw}
T_1(z)=e^{-|x|}{\Ke}(z, z)=2^{n-1}e^{-|x|}
\sum_{j=0}^{\infty}\frac{(j+n-2)!}{(2j+n-2)!}\frac{|x|^{2j}}{j!},
\ 2|x|^2=z\cdot \bar z .
\end{equation}
Now, we compute the Kempf distortion functions $T_{m}(z)$ for all
non-negative integers
$m$ as follows.
Making the change of variable $mz=w$, we get
$$\|f\|_m^2 = \int_X |f(w/m)|^2 e^{-| \PI w|}
m^{-n} d\mu (w),$$
since $d\mu(w/m)=m^{-n}d\mu(w)$.
Consequently, the operator
$$ T:\ Tf(w):=m^{-\frac{n}{2}} f(w/m)  $$
is a unitary isomorphism from ${\cal H}_m$ onto ${\cal H}$
\footnote{The second author is in debt with
Miroslav Engli\v{s} who pointed him out
the idea of using this isomorphism to  compute
$T_{m}$ from  $T_1$.}.
Denoting by
$\Ke^{(m)}(w,z)\equiv\Ke^{{m}}_z(w)$ the reproducing kernel of ${\cal H}_m$
(and writing simply $\Ke (w,z)\equiv K_z(w)$ if $m=1$), we therefore have,
on the one hand,
$$f(z) = \langle f,{\Ke} ^{(m)}_z\rangle_m = \langle Tf, T{\Ke}
^{(m)}_z\rangle  $$
for any $f\in {\cal H}_m$, while, on the other hand,
$$ f(z) = m^{\frac{n}{2}} Tf(mz) = \langle Tf,  m^{\frac{n}{2}}
{\Ke}_{mz}\rangle  .$$
Thus $T\Ke^{(m)}_z= m^{\frac{n}{2}} {\Ke}_{mz}$, and
$${\Ke}^{(m)}_z(w) =  m^{\frac{n}{2}} T^{-1} {\Ke}_{mz}(w) = m^{n}
{\Ke}_{mz}(mw).  $$
That is,
$${\Ke}^{(m)}(w,z) = m^{n} {\Ke}(mw,mz) .  $$

Substituting this into Rawsley's formula
(\ref{raw}), we thus get

\begin{eqnarray}\label{Tm}
T_m(z)& = & e^{-m|x|}{\Ke}^{(m)}(z,z) =
2^{n-1}m^ne^{-m|x|}\sum_{j=0}^{\infty}\frac{(j+n-2)!}{(2j+n-2)!}\frac{(m|x|)
^{2j}}{j!} 
\end{eqnarray}

\begin{remar}\rm
From (\ref{Tm}) one sees that
$T_{m}(x) = m^{n} T_{1}(mx)$.
In the compact case the
relationship between $T_{m}$
and $T_1$ is unknown.
Here the fact that the Hilbert spaces involved
are infinite dimensional is a crucial step
to get the previous equality.
\end{remar}

Notice that from representation (\ref{Tm}) is not clear the growth of $T_m(z)$
as $m\rightarrow\infty$.
The following proposition gives us important analytic information about $T_m$
as $m\rightarrow\infty$.

\begin{proposition}
Kempf's distortion function for the Kepler manifold
can be written in the following two forms: 
\begin{equation}\label{Tm2}
T_m(z) = 
2^{-1} m^ne^{-m|x|}\sum_{j=0}^{\infty} (1 + \tau_j)  
\frac{(m|x|)
^{2j}}{(2j)!},
\end{equation}
where
$$\tau_j =  
1- \frac{(j+1)\ldots (j+n-2)}{j+1/2)\ldots (j+(n-2)/2} \longrightarrow 0 \ \ \ \textrm{ for $j\rightarrow \infty$,}$$
and 
\begin{equation}\label{Tmnew}
T_m(z)=
2m^ne^{-\xi_m}
(\frac{1}{\xi_m}\frac{\partial}{\partial \xi_m})^{n-2}
[\xi_m^{n-2}(e^{\xi_m}+(-1)^{n-2}e^{-\xi_m}+ Q(\xi_m))],
\end{equation}
where $\xi_m=m|x|$, $Q(\xi_m)$ is a polynomial of degree $\leq n-4$ in the
variable $\xi_m$.
\end{proposition} 
\textit{Proof}.
From (\ref{Tm}) one gets
\begin{eqnarray}
T_m(z)& = & e^{-m|x|}{\Ke}^{(m)}(z,z) =
2^{n-1}m^ne^{-m|x|}\sum_{j=0}^{\infty}\frac{(j+n-2)!}{(2j+n-2)!}\frac{(m|x|)
^{2j}}{j!}\nonumber\\
& = & 
2^{n-1}m^ne^{-m|x|}\sum_{j=0}^{\infty}\frac{(j+n-2)! (2j)!}{j!(2j+n-2)!}\frac{(m|x|)
^{2j}}{(2j)!}\nonumber\\
& = & 
2^{-1} m^ne^{-m|x|}\sum_{j=0}^{\infty}\frac{(j+1)\ldots (j+n-2)}{j+1/2)\ldots (j+(n-2)/2} 
\frac{(m|x|)
^{2j}}{(2j)!}\nonumber\\ 
& = & 
2^{-1} m^ne^{-m|x|}\sum_{j=0}^{\infty} (1 + \tau_j)  
\frac{(m|x|)
^{2j}}{(2j)!},
\end{eqnarray}
In order to prove (\ref{Tmnew}) set
$$y_m^2=m|x|=\xi_m, \ y_m, \xi_m\in {\real}\setminus\{0\}.$$
Then, since $\frac{\partial}{\partial y_m}=
\frac{1}{2\xi_m}\frac{\partial}{\partial \xi_m}$ one gets 
$$\begin{array}{lll}T_{m}(z)&=&
2^{n-1}m^ne^{-\xi_m}\sum_{j=0}^{\infty}\frac{(j+n-2)!}{(2j+n-2)!}\frac{y_m^{
j}}{j!}\\
&=& 2^{n-1}m^ne^{-\xi_m}(\frac{\partial}{\partial y_m})^{n-2}
\sum_{j=0}^{\infty}\frac{y_m^{j+n-2}}{(2j+n-2)!}\\
&=&2 m^ne^{-\xi_m}(\frac{1}{\xi_m}\frac{\partial}{\partial \xi_m})^{n-2}
[\xi_m^{n-2}\sum_{j=0}^{\infty}\frac{\xi_m^{2j+n-2}}{(2j+n-2)!}]
\end{array}$$

If $n$ is even then
$$T_{m }(z)=
2m^ne^{-\xi_m}
(\frac{1}{\xi_m}\frac{\partial}{\partial \xi_m})^{n-2}
[\xi_m^{n-2}(\cosh \xi_m- P(\xi_m))],$$
where
$$P(\xi_m)=\sum_{j=0}^{\frac{n-4}{2}}\frac{\xi_m^{2j}}{(2j)!}.$$
If $n$ is odd then
$$T_m(z)=
2m^ne^{-\xi_m}
(\frac{1}{\xi_m}\frac{\partial}{\partial \xi_m})^{n-2}
[\xi_m^{n-2}(\sinh \xi_m- R(\xi_m))],$$
where
$$R(\xi_m)=\sum_{j=0}^{\frac{n-5}{2}}\frac{\xi_m^{2j+1}}{(2j+1)!},$$
and hence (\ref{Tmnew}) easily follows.
\hfill $\Box$

\begin{remar}\rm
By (\ref{Tm2}) one might view $T_m(z)$ as a ``small perturbation'' (or close to) for $m\rightarrow \infty$ of 
$$T_m^0(z)=  
2^{-1} m^ne^{-m|x|}\sum_{j=0}^{\infty}   
\frac{(m|x|)
^{2j}}{(2j)!} = m^ne^{-m|x|} \cosh (m|x|) = m^n \frac{1-e^{-2m|x|}}{2}.$$
\end{remar}

\section{TYZ expansion for the Kepler manifold}\label{expansion}
The key ingredient 
to find the TYZ expansion 
of $T_m$ for the Kempf distortion 
function of the Kepler manifold
is (\ref{Tmnew}). 
Clearly we have 
\begin{eqnarray}
T_m(z)  =  
2m^n F (m|x|),  \label{Trepr1}
\end{eqnarray}
where
\begin{eqnarray}
F (y)  = 
e^{-y}
(\frac{1}{y}\frac{d}{d y } )^{n-2}
\left( y^{n-2}(e^{y}+(-1)^{n-2}e^{-y}+ Q(y))\right),   \qquad y\in \R .
\label{Trepr2}
\end{eqnarray}

The explicit representation \eqref{Trepr1}-\eqref{Trepr2} of  $T_m(z)$ for the Kepler manifold 
has a remarkable feature, namely, it is defined by a generating function $F(y)$ depending on one variable. 
Note that in fact $T_m(z)$ is independent of the base variables $e \in S^n$.

The first main result of the present paper is the following one.

\begin{theorem}\label{mainteor}
Let $F$ satisfy \eqref{Trepr2}. Then the following representation holds: 
\begin{eqnarray}
F(y) &= & \sum_{j=0}^{n-2}  \frac{b_j}{y^j}  + \Phi(y) + \Psi(y)
\label{rep1}
\end{eqnarray}
where
\begin{eqnarray}
\Phi(y) 
& = & e^{-2y} \sum_{j=0}^{n-2}  \frac{p_j}{y^j} \label{rep2} \\
\Psi(y) & = & e^{-y} \sum_{j=0}^{n-3}  \frac{r_j}{y^j} 
\label{rep3}
\end{eqnarray}
and the constants $a_j$, $p_j$, $r_j$ are calculated explicitly. 
The functions $\Phi(y)$, $\Psi(y)$ and therefore,  $F(y)$ as well,  are extended 
to holomorphic functions in semiplane $\Re y >0$. 
In particular, by  \eqref{Trepr1} and \eqref{rep1}
we get 
\begin{eqnarray}
T_m(z) &= & \sum_{j=0}^{n-2}  a_j(x) m^{n-j}  + 2m^n\Phi(m|x|) + 2m^n\Psi(m|x|),
m\in \N,
\label{rep1T}
\end{eqnarray}
where
\begin{eqnarray}
a_j(x) &= &   \frac{2b_j}{|x|^j}, \qquad j=0,1,\ldots, n-2.  
\label{rep1Ta}
\end{eqnarray}
and
\begin{eqnarray}
a_0(x)& = & 1
\label{expcoef0} \\
a_1(x) & = & \frac{(n-2)(n-1)}{2|x|} \label{expcoef1}
\label{expcoef2}
\end{eqnarray}
Moreover, 
there exists an absolute constant $C_0>0$ such that for every $\delta \in ]0,1]$ 
\begin{eqnarray}
\sup_{|x| \geq \delta }|D_x^\alpha \Theta_m(x)|,  & \leq & C_0^{\alpha+1} \frac{\alpha!}{\delta^\alpha } 
e^{-m\delta/2}
\label{TYZ1}
\end{eqnarray}
for all $m\in \N$, where $\Theta =\Phi, \Psi$. Therefore, we have the following estimates
\begin{eqnarray}
|D_x^\alpha \left( T_{m} - \sum_{j=0}^{n-2} a_j(x) m^{n-j} \right)| & \leq & 
C_0^{\alpha+1}  \frac{\alpha!}{\delta^\alpha }e^{-m\delta/2}
\label{TYZ2}
\end{eqnarray}
for all $|x| \geq \delta$, $\alpha \in \Z_+^n$. 
\end{theorem}

\textit{Proof}. We recall the well known Fa\`a di Bruno  type formula for the derivative of $g\circ \varphi$,  
namely, for a given $\alpha\in\N$ we have
\begin{eqnarray}
\nonumber
D_t^\alpha (g (\varphi (t))) & = &  D_y^\alpha (g(\varphi (y)))|_{y=t}\\
&=&
\label{faa1}
 \sum_{j=1}^\alpha \frac{g^{(j)}(\varphi(t))}{j!} D_x^\alpha \left( (\varphi(x)-\varphi(t))^j\right)|_{x=t}\\
& = & \sum_{j=1}^\alpha \frac{g^{(j)}(\varphi(t))}{j!}  
\sum_{\stackrel{\alpha_1+\cdots+\alpha_j=\alpha}{\al_1\geq 1, \ldots,\alpha_j\geq 1}}
\frac{\alpha!}{\alpha_1!\ldots \alpha_j!} \varphi^{(\alpha_1)}(t) \ldots \varphi^{(\alpha_j)}(t), 
\nonumber
\end{eqnarray}
where $\varphi^{(k)}(t)$ stands for $D_t^k\varphi (t)$.

Next, we straighten $y^{-1}D_y$ into $D_t$ via  the change of the variable $y =y(t) = \sqrt{2t}$, $t=t(y) = y^2/2$. 
Therefore, setting 
\begin{equation}
G(t) = F(\sqrt{2t}), \ t>0, \qquad F(y) = G(\frac{y^2}{2}), \ y>0, 
\label{pt2} 
\end{equation}
we get by 
\eqref{Trepr2}
\begin{equation}
F(y) = G(t) = e^{-\sqrt{2t} } \left( \frac{d}{dt}\right)^{n-2} \left( (2t)^{(n-2)/2} \frac{e^{\sqrt{2t} }+ (-1)^{n-2}e^{-\sqrt{2t}}}{2} + Q(\sqrt{2t})\right).
\label{exloi1}
\end{equation}

The next assertion is instrumental in the proof.

\begin{lemma}
Let $N\in \N$, $c\in \R$, and $r>0$. Then 
\begin{eqnarray}
\psi_N^{c,r}(y) &:=& e^{-y} \left( \frac{1}{y}\frac{d}{dy}\right)^N ( y^r e^{cy} ) \nonumber\\
& = & \psi_N^{c,r}(\sqrt{2t}) = :\varphi_N^{c,r}(t)= 
e^{-\sqrt{2t} } \left( \frac{d}{dt}\right)^N ( (2t)^{r/2} e^{c\sqrt{2t} } ).
\label{exloi2}
\end{eqnarray}
has the following representation
\begin{eqnarray}
\varphi_N^{c,r}(t)& = & 
e^{-(1-c)\sqrt{2t} } (2t)^{(r-N)/2} \sum_{s=0}^{N} \frac{\varkappa_s}{(2t)^{s/2}},  
\label{exloi3t}
\end{eqnarray}
i.e. 
\begin{eqnarray}
\psi_N^{c,r}(t)& = & 
e^{-(1-c)z} z^{(r-N)/2} \sum_{s=0}^{N} \frac{\varkappa_s}{z^s},  
\label{exloi3z}
\end{eqnarray}
where 
\begin{eqnarray}
\varkappa_s& = & 
\frac{1}{ (N-s)! } 
 \sum_{\ell =N-s}^N  \left( \begin{array}{c}N\\ \ell \end{array}\right)  
 \left( \prod_{q=0}^{N-\ell -1} (\frac{r}{2}-q) \right)  2^{N-r/2} (-1)^{\ell+s-N}  \nonumber\\
 & \times &  
\sum_{\stackrel{\ell_1+\cdots+\ell_{N-s}=\ell}{\ell_1\geq 1, \ldots,\ell_{N-s}\geq 1}}\, 
\frac{\ell!}{\ell_1!\ldots \ell_{N-s}!} \prod_{q_1}^{\ell_1-1} (\frac{1}{2}-q_1)  \ldots \prod_{q_{N-s}}^{\ell_{N-s}-1} (\frac{1}{2}-q_{N-s})
\label{exloi4}
\end{eqnarray}
for $s=0,\ldots, N-1$ and 
\begin{eqnarray}
\varkappa_N& = & 2^{N-r/2} \prod_{q=0}^{N-1}
(\frac{r}{2}-q).
\label{exloi4a}
\end{eqnarray}
\end{lemma}

\textit{Proof}. By  Fa\`a di Bruno type formula \eqref{faa1} we derive
\begin{eqnarray}
\Theta^{r,c}_N(t) &=& \left( \frac{d}{dt}\right)^N ( t^{r/2} e^{c\sqrt{2t} } )\nonumber\\
 & =& \sum_{\ell=0}{N} \left( \begin{array}{c}N\\ \ell \end{array}\right) D_t^{N-\ell}(t^{r/2}) 
 D_t^\ell (e^{c\sqrt{2t} })\nonumber\\ 
 & = &   D_t^N(t^{r/2}) e^{c\sqrt{2t} } +  
 \sum_{\ell=1}N \left( \begin{array}{c}N\\ \ell \end{array}\right)  
 \left( \prod_{q=0}^{N-\ell -1} (\frac{r}{2}-q) \right)  t^{r/2-N+\ell}
 e^{c\sqrt{2t}}\sum_{j=1}^\ell \frac{(c2)^{j/2}}{j!}  \nonumber\\
 & \times &  
\sum_{\stackrel{\ell_1+\cdots+\ell_j=\ell}{\ell_1 \geq 1, \ldots, \ell_j \geq 1} }\, 
\frac{\ell!}{\ell_1!\ldots \ell_j!} D_t^{\ell_1}(t^{1/2}) \ldots D_t^{\ell_j} (t^{1/2})
\label{pfaa1}
\end{eqnarray}
with the convention $\prod_{q=0}^{-1} ... = 1$. 
Since 
\begin{eqnarray}
D_t^\mu (t^{1/2}) & = & \frac{1}{2} (\frac{1}{2}-1) \ldots (\frac{1}{2} -\mu+1) t^{1/2-\mu}\nonumber\\ 
& = & (-1)^{\mu-1}\frac{(2\mu-3)!!}{2^\mu} t^{ 1/2 -\mu }
\label{pfaa2}
\end{eqnarray}
for all positive integers $\mu$, with $(-1)!!:= 1$, $(2\mu-3)!! := 1 \ldots (2\mu-3)$ if $\mu \geq 2$, 
combining \eqref{pfaa1} and \eqref{pfaa2}, we obtain
\begin{equation}
\sum_{\stackrel{\ell_1+\cdots+\ell_j=\ell}{\ell_1\geq 1, \ldots,\ell_j\geq 1}}\, 
\frac{\ell!}{\ell_1!\ldots \ell_j!} D_t^{\ell_1}(t^{1/2}) \ldots D_t^{\ell_j} (t^{1/2}) = 
(-1)^{\ell-j} \Gamma^{\ell,j} 
\frac{(2t)^{j/2 -\ell} }{2^{j/2}}  
\label{pfaa3}
\end{equation}
with 
\begin{equation}
\Gamma^{\ell,j} := 
\sum_{\stackrel{\ell_1+\cdots+\ell_j=\ell}{\ell_1\geq 1, \ldots,\ell_j\geq 1}}\, 
\frac{\ell!}{\ell_1!\ldots \ell_j!} (2\ell_1-3)!! \ldots (2\ell_j-3)!!.
\label{pfaa4}
\end{equation}
We note that 
\begin{eqnarray}
\Gamma^{\ell,\ell} & = & \ell! 
\label{pfaa41}\\
\Gamma^{\ell,\ell-1} &  = &  -\frac{\ell-1}{2} \ell! 
\label{pfaa42}
\end{eqnarray}

Therefore, by \eqref{pfaa1} - \eqref{pfaa3},  
\begin{eqnarray}
\Theta^{r,c}_N(t)& = &  2^{N-r/2} \left( \prod_{q=0}^{N-1}
(\frac{r}{2}-q) \right) (2t)^{r/2-N}
e^{c\sqrt{2t} } \nonumber\\
& + & 
 \sum_{\ell=1}^N \left( \begin{array}{c}N\\ \ell \end{array}\right)  
 \left( \prod_{q=0}^{N-\ell -1} (\frac{r}{2}-q) \right)  2^{N-\ell -r/2} 
 (2t)^{r/2-N+\ell} \nonumber\\
 & \times &  \sum_{j=1}^\ell \frac{c^{j/2}}{j!}(-1)^{\ell-j} 
 \Gamma^{\ell,j} (2t)^{j/2 -\ell} \nonumber\\
 & = & (2t)^{r/2-N/2}
e^{c\sqrt{2t} }
\frac{ 2^{N-r/2} \prod_{q=0}^{N-1}
(\frac{r}{2}-q) }{ (2t)^{N/2} }
 \nonumber\\
& + & (2t)^{r/2-N/2}e^{c\sqrt{2t} } 
 \sum_{j=1}^N  \frac{1}{ j!(2t)^{(N -j)/2} } 
 \sum_{\ell =j}^N  \left( \begin{array}{c}N\\ \ell \end{array}\right)  
 \left( \prod_{q=0}^{N-\ell -1} (\frac{r}{2}-q) \right)  2^{N-\ell -r/2} (-1)^{\ell-j} 
 \Gamma^{\ell,j}\nonumber\\
 & = & (2t)^{r/2-N/2}
e^{c\sqrt{2t} }
\frac{ 2^{N-r/2} \prod_{q=0}^{N-1}
(\frac{r}{2}-q) }{ (2t)^{N/2} }
 \nonumber\\
& + & (2t)^{r/2-N/2}e^{c\sqrt{2t} } 
 \sum_{s=0}^{N-1}  \frac{1}{ (N-s)!(2t)^{s/2} } 
 \sum_{\ell =N-s}^N  \left( \begin{array}{c}N\\ \ell \end{array}\right) \nonumber\\
 & \times &  
 \left( \prod_{q=0}^{N-\ell -1} (\frac{r}{2}-q) \right)  2^{N-\ell -r/2} (-1)^{\ell+s-N} 
 \Gamma^{\ell,N-s}\nonumber\\
 & = & (2t)^{r/2-N/2} e^{c\sqrt{2t} } \sum_{s=0}^{N} \frac{\varkappa_s}{(2t)^{s/2}}
 \label{pfaa5}
\end{eqnarray}
where $\varkappa_s$ is defined by 
\begin{eqnarray}
\varkappa_s & := & \frac{1}{ (N-s)! } 
 \sum_{\ell =N-s}^N  \left( \begin{array}{c}N\\ \ell \end{array}\right)  
 \left( \prod_{q=0}^{N-\ell -1} (\frac{r}{2}-q) \right)  2^{N-\ell -r/2} (-1)^{\ell+s-N}  
 \Gamma^{\ell,N-s}\nonumber
\end{eqnarray}

In view of the definition of $\Gamma^{\ell,j}$ with the convention $\Gamma^{\ell,0}=1$, it is equivalent to \eqref{exloi4},  \eqref{exloi4a}. This ends the proof of the lemma. \hfill $\Box$
\vspace{10pt}

We conclude the proof of the theorem by applying the
previous  lemma  for $z= m|x|$ and obtain the value of  $a_s = 1/2 \varkappa_s^{N,r;c}$ by setting  
$c=1$, $r=N=(n-2)$; $p_s= (-1)^{n-2}/2 \varkappa_s^{N,r;c}$ by setting $c=-1$, $r=N=n-2$; and 
$$ r_s = \sum_{j=0}^{n-3} q_j\varkappa^{n-2, j;0} $$ provided $\ds Q(z) = \sum_{j=0}^{n-3} q_jz^j$. 
 \hfill $\Box$.

\vspace{10pt}

\begin{remark}
In view of \eqref{rep1T}, we have 
\begin{eqnarray}
T_m(z) &= &  m^n  + \frac{(n-2)(n-1)}{2 |x|} m^{n-1} +
 \sum_{k=2}^{n-2}  \frac{2a_k}{|x|^k} m^{n-k}  + R_m(|x|),
\label{rep1Ta1}
\end{eqnarray}
with $R_m(x) $ being exponentially small $e^{-cm}$ away from the origin $x=0$.  
\end{remark}

\begin{remark}
The novelty of the theorem above is twofolded.
First, our TYZ type expansion is finite, i.e.,
$a_j=0$ for $j\geq n-1$ (compare (\ref{rest})).
Secondly, the reminder is exponentially small.
Moreover,  the coefficients $a_j$
can be computed explicitely.
In a forthcoming paper we study tha link between these coefficients
$a_j$ and the curvature of the metric $g$ as in Lu's Theorem
\cite{lu}.
\end{remark}

\begin{remar}
One can also investigate the asymptotic expansion near the singular (conic)
point.
Using a local coordinates in which it coincides   with the origin, we can
derive explicit  asymptotic expansion for Kempf's distortion function
$T_m  ( x )$ near $x=0$. Moreover one can show that
\begin{equation}
\| T_m  ( x ) - \sum_{j=0}^{n-2} a_j(x) m^{n-j}\|_{L^p( B(\delta))} =O(
\delta^{n-(n-2)p}) ) m^2, \ \ \ \delta \searrow 0, m\geq 1
\label{TYZsing1}
\end{equation}
provided $1\leq p < n/(n-2)$. So we encounter the critical $L^p$ index which
appears
in different mathematical problems.
\end{remar}

\section{Proof that our estimate is sharp}\label{sharp}
As a consequence of Theorem \ref{mainteor} and 
\ref{obstrer1} the K\"{a}hler form
$g$ on the Kepler manifold $X$
is the $C^{\infty}$-limit of suitable
normalized projectively induced K\"{a}hler
metrics,
namely
$$\lim_{m\rightarrow\infty}\frac{1}{m}\varphi_m^*(g_{FS})=g$$
where
$\varphi_m: X\rightarrow {\complex}P^{\infty}$
is the coherent states map.
In this  Section
we show that $g$ is not projectively induced
(via any map) and  then that our extimate 
in  Theorem \ref{mainteor} is sharp.

We need to recall briefly 
some results
about Calabi's diastasis function  referring  the reader to
\cite{ca} and \cite{diastherm} for details and further results.

Let $M$ be a complex manifold
endowed with a
real analytic K\"{a}hler metric $g$.
Then, in a neighborhood of every point
$p\in M$, one can introduce
a very special
K\"{a}hler potential
$D^g_p$ for the
K\"{a}hler form
$\omega$ associated to $g$, which
Calabi \cite{ca} christened
{\em diastasis}.
Recall that
a K\"{a}hler potential
is an
analytic function
$\Phi$
defined in a neighborhood
of a point $p$
such that
$\omega =\frac{i}{2}\bar\partial\partial\Phi$.
A K\"{a}hler potential is not unique:
it is defined up to an addition with
the real part of a holomorphic function.
By duplicating the variables $z$ and $\bar z$
a potential $\Phi$ can be complex analytically
continued to a function
$\tilde\Phi$ defined in a neighborhood
$U$ of the diagonal containing
$(p, \bar p)\in M\times\bar M$
(here $\bar M$ denotes the manifold
conjugated to $M$).
The {\em diastasis function} is the
K\"{a}hler potential $D^g_p$
around $p$ defined by
$$D^g_p(q)=\tilde\Phi (q, \bar q)+
\tilde\Phi (p, \bar p)-\tilde\Phi (p, \bar q)-
\tilde\Phi (q, \bar p).$$

Observe that the diastasis does not depend
on the potential chosen,  $D^g_p(q)$
is symmetric in $p$ and $q$
and $D^g_p(p)=0$.

The  diastasis function is the key
tool for studying
the K\"{a}hler
immersions of a K\"{a}hler manifold into
another K\"{a}hler manifold
as expressed by the following lemma.

\begin{lemma}(Calabi \cite{ca})\label{dxy}
Let $(M, g)$ be a K\"{a}hler manifold
which admits a K\"{a}hler immersion
$\varphi:(M, g)\rightarrow (S, G)$ into a real
analytic K\"{a}hler manifold
$(S, G)$. Then $g$ is real analytic.
Let $D^g_p:U\rightarrow {\real}$
and $D^G_{\varphi(p)}:V\rightarrow {\real}$ be
the diastasis functions of $(M, g)$ and $(S, G)$
around $p$ and $\varphi(p)$ respectively.
Then  $\varphi^{-1}(D^{G}_{\varphi(p)})=D^{g}_p$
on $\varphi^{-1}(V)\cap U$.
\end{lemma}

When $(S, G)$
is the $N$-dimensional complex projective space
$S={\complex}P^N$
equipped with
with the Fubini--Study metric $G=g_{FS}$,
one can show that
for all $p\in {\complex}P^N$
the diastasis function $D^{g_{FS}}_{p}$
around $p$ is globally defined
except in the cut locus
$H_p$ of $p$ where it blows up.
Moreover, $e^{-D^{g_{FS}}_{p}}$
is globally defined  (and smooth) on
${\complex}P^N$
(see \cite{ca} or \cite{diastherm} for details).

Then, 
by Lemma \ref{dxy} one immediately  gets the following:

\begin{lemma}\label{dxyz}
Let  $g$ be a projectively induced K\"{a}hler metric
on a  complex manifold $M$.
Then, $e^{-D^{g}_p}$ is globally defined on all $M$.
\end{lemma}

\begin{corol}\label{coroln}
Let $g_*$ be the  K\"{a}hler metric on
${\complex}^*$ whose associated K\"{a}hler form is given by
$\omega_*=\frac{i}{2}\partial\bar\partial |\eta|, \eta=x+iy$.
Then $g_*$ is not projectively induced.
\end{corol}
\dimostr
Fix any point $\alpha\in {\complex}^*$.
A globally defined K\"{a}hler potential $\Phi$
for the K\"{a}hler metric $g_*$
around  $\alpha$
is given by $\Phi (\eta)=|\eta|$
and Calabi's diastasis function
around $\alpha$ reads as
$$D^{g_*}_{\alpha}:U\rightarrow{\real},\
\eta\mapsto |\eta|+|\alpha|-\sqrt{\eta\bar \alpha}-\sqrt{\bar \eta
\alpha},$$
where $U\subset{\complex}^*$ is suitable
simply-connected open subset of
${\complex}^*$ around $\alpha$ (as a maximal domain
of defintion of  $D^{g_*}_{\alpha}$
one can take $U={\complex}^*\setminus L$
where $L$ is any half-line starting from the origin of
${\complex}={\real}^2$
such that $\alpha\notin L$).
The function  $D^{g_*}_{\alpha}$,
as well as  the function
$e^{-D^{g_*}_{\alpha}}$,
cannot be extended to all ${\complex}^*$.
Hence  we are done
by Lemma \ref{dxyz}.
\fdim

We are now in the position to prove
that our estimate  is sharp.

\begin{teor}\label{mainteor2}
Let $g$ be the K\"{a}hler metric on the Kepler
manifold $X$
whose associated K\"{a}hler form is given by
(\ref{Kform}).
Then $g$ is not projectively induced.

\end{teor}
\dimostr
First observe that the map
$$j :({\complex}^*, g_* )
\rightarrow (X, g)$$
defined by
$j (z)=
(\eta, i\eta, 0, \dots, 0 )$
is a K\"{a}hler immersion
satisfying
$j ^{*}(g)=g_*$, 
with $g_*$ as in Corollary \ref{coroln}.
Assume by contradiction   that
$g$ is projectively induced,
namely  there exists
$N\leq\infty$
and a  K\"{a}hler
immersion
$\varphi:(X, g )\rightarrow ({\complex}P^N, g_{FS})$.
Then the map
$\varphi\circ j:({\complex}^*, g_*)\rightarrow ({\complex}P^N, g_{FS})$
would be a  K\"{a}hler immersion
contradicting  Corollary \ref{coroln}.
\fdim




\section{Estimates of the logarithmic obstruction term}\label{sharplog}

The aim of this section is two--folded. First, taking advantage of the 
homogeneity  structure  of the Kepler manifold 
we introduce global polar--angular coordinates. As an outcome, 
we are able to write down  the explicit form of 
the operator   $ \partial\bar\partial $
 in such coordinates. 
Secondly, we show  precise asymptotic expansion and estimates of the logarithmic error term. 
 This is a novel result, as far as we know. 
A key functional--analytic ingredient of our arguments in the proof of the representation formula in the homogeneous coordinates  
is the fact that the distortion function for 
the Kepler manifold $X$ depends only on the rescaled by the factor $\sqrt{2}$ Euclidean distance to the origin in $C^{n+1}$ identified with $\R^{2n+2}$ by the canonical complex structure $J_0$.  
More precisely, the Kepler manifold is given by the isotropic cone $C$, defined globally by 
\begin{eqnarray}
C & =  &]0, +\infty[ \times \widetilde{C}, \label{polar1}
\end{eqnarray}
where
\begin{eqnarray}
\widetilde{C} & = & 
\{ e +is \in \C^{n+1}: \,    (e, s) \in S^n\times S^n:\, \ e\cdot s =0\}
\label{polar2}
\end{eqnarray} 
is the unitary tangent bundle of $S^n$.  Here 
$e= (e_1,\ldots, e_{n+1})$, $s= (s_1,\ldots, s_{n+1})$, and 
$e\cdot s=\sum_{j=1}^{n+1}e_js_j$. 
We can parameterize explicitly $\widetilde{C} $ by using twice the classical angular variables 
in $\R^{n+1}$, namely, $e =e (\varphi)$ 
and  $s =s(\psi)$, with  
\begin{eqnarray}
e_1 & =  & \sin \varphi_1 \ldots \sin  \varphi_{n-1}  \sin \varphi_n,  \nonumber\\
e_2 & =  & \sin \varphi_1 \ldots \sin  \varphi_{n-1}  \cos  \varphi_n,  \nonumber\\
e_3 & =  & \sin \varphi_1 \ldots  \sin\varphi_{n-2} \cos  \varphi_{n-1},  \nonumber\\
\ldots & =& \ldots\nonumber\\
e_{n+1} & = & \cos \varphi_1,
\label{polar1n}
\end{eqnarray}
and 
\begin{eqnarray}
s_1 & =  & \sin \psi_1 \ldots \sin  \psi_{n-1}  \sin \psi_n,  \nonumber \\
s_2 & =  & \sin \psi_1 \ldots \sin  \psi_{n-1}  \cos  \psi_n,  \nonumber \\
s_3 & =  & \sin \psi_1 \ldots  \sin\psi_{n-2} \cos  \psi_{n-1},  \nonumber \\
\ldots & =& \ldots\nonumber \\
s_{n+1} & = & \cos \psi_1, 
\label{polar2n}
\end{eqnarray}
with $\varphi_1 \in 
[0, + 2\pi[$, $\psi_1 \in
[0, 2\pi[$, $\varphi_j\in [0, \pi[$, 
$\psi_j \in [0,  \pi[$, 
$j=2,\ldots, n$.

Clearly  $\widetilde{C} \subset S^{2n+1}(\sqrt{2})$ since 
\begin{equation}
e+is  \in  S^{2n+1}(\sqrt{2}) = \{ \zeta \in \C^{n+1}: \, \zeta , \bar\zeta = 2 \} , 
\qquad \textrm{if  $e+is\in \widetilde{C}$}.
\label{emb1}
\end{equation}

If we define $S^{2n+1}(r)$ by the angular coordinates in $\R^{2n+2}=\C^{2n+2}$ using the coordinates $\zeta = \xi + i\eta \in \C^{n+1}$ and 
$(\xi_1,\eta_1, \ldots, \xi_{n+1}, \eta_{n+1})\in \R^{2n+2}\setminus 0$ (identifying with 
$\zeta = \xi + i\eta \in \C^{n+1}\setminus 0$) by the standard angular cordinates  
we will have that $S^{2n+1}(\sqrt{2})$ is defined by $r=\sqrt{2}$, where 
\begin{eqnarray}
\xi_1 & =  & r\sin \theta_1 \ldots \sin  \theta_{2n}  \sin \theta_{2n+1}  \nonumber\\
\eta_1 & =  & r\sin \theta_1 \ldots \sin  \theta_{2n}  \cos  \theta_{2n+1}  \nonumber\\
\ldots & =& \ldots\nonumber\\
\xi_{n+1} & = & r\sin \theta_1 \cos \theta_2 \nonumber\\
\eta_{n+1} & = & r\cos \theta_1
\label{polar1can}
\end{eqnarray}
with $r>0$, $\theta_1 \in 
[0, + 2\pi[$, $\theta_j\in [0, \pi[$,  
$j=2,\ldots, 2n+1$. However, such coordinates do not provide an easy definition of  $\widetilde{C}$ by implicit function theorem. 

We construct an embedding  $\widetilde{C} $ in $S^{2n+1}(\sqrt{2})$ 
{\em compatible} with the standard complex structure $J_0$ of  
 $\C^{n+1}$ by introducing apparently new angular coordinates  on $S^{2n+1}(\sqrt{2})$ which differ from the traditional ones. 
More precisely, we set
\begin{eqnarray}
\xi_1 & =  & r\sin \varphi_1 \ldots \sin  \varphi_{n-1}  \sin \varphi_n\cos \theta  \nonumber\\
\xi_2& =  & r\sin \varphi_1 \ldots \sin  \varphi_{n-1}  \cos  \varphi_n \cos \theta   \nonumber\\
\xi_3 & =  & r\sin \varphi_1 \ldots  \sin\varphi_{n-2} \cos  \varphi_{n-1} \cos \theta \nonumber\\
\ldots & =& \ldots\nonumber\\
\xi_{n+1} & = & r\cos \varphi_1 \cos \theta, 
\label{angemb1a}
\end{eqnarray}
\begin{eqnarray}
\eta_1 & =  & r\sin \psi_1 \ldots \sin  \psi_{n-1}  \sin \psi_n \sin \theta  \nonumber \\
\eta_2 & =  & r\sin \psi_1 \ldots \sin  \psi_{n-1}  \cos  \psi_n \sin \theta \nonumber \\
\eta_3 & =  & r\sin \psi_1 \ldots  \sin\psi_{n-2} \cos  \psi_{n-1}\sin \theta  \nonumber \\
\ldots & =& \ldots\nonumber \\
\eta_{n+1} & = & r\cos \psi_1 \sin \theta,  
\label{angemb1b}
\end{eqnarray}
where $r>0$, and $\varphi = (\varphi_1, \ldots, \varphi_n) \in I$, 
$\psi = (\psi_1,\ldots, \psi_n)\in J$, $\theta \in [0,2\pi[$,   
\begin{eqnarray}
I & =  & I_1\times I_2\times \ldots I_n \label{defpol1}\\
J & =  & J_1\times J_2\times \ldots J_n \label{defpol2}
\end{eqnarray}
where $I_1,J_1$ are semi--closed intervals of length $2\pi$ while 
$I_2,\ldots,I_n, J_2, \ldots, J_n$ 
are semi--closed intervals of length $\pi$.

\vspace{10pt}

\begin{example}
Let $n=1$. Then the polar coordinates \eqref{angemb1a}, \eqref{angemb1b} for $\C^2 = \R^4$ 
become
\begin{eqnarray}
\xi_1 & =  & r\sin \varphi_1 \cos \theta,\nonumber\\
\xi_2 & = & r \cos \varphi_1 \cos \theta 
\label{angemb1a-d2}
\end{eqnarray}
and 
\begin{eqnarray}
\eta_1 & =  & r\sin \psi_1 \sin \theta, \nonumber \\
\eta_2 & = & r\cos \psi_1 \sin \theta. 
\label{angemb1b-d2}
\end{eqnarray}
\vspace{10pt}

If  $n=3$, the polar coordinates \eqref{angemb1a}, \eqref{angemb1b} for $\C^3 = \R^6$ 
become 
\begin{eqnarray}
\xi_1 & =  & r\sin \varphi_1  \sin \varphi_2 \cos \theta  \nonumber\\
\xi_2& =  & r\sin \varphi_1  \cos \varphi_2  \cos \theta   \nonumber\\
\xi_3 & =  & r\sin \cos\varphi_1  \cos \theta \nonumber\\
\label{angemb1a-d3}
\end{eqnarray}
and 
\begin{eqnarray}
\eta_1 & =  & r\sin \psi_1 \sin  \psi_2  \sin \theta  \nonumber \\
\eta_2 & =  & r\sin \psi_1 \cos  \psi_2 \sin \theta \nonumber \\
\eta_3 & = & r\cos \psi_1 \sin \theta, 
\label{angemb1b-d3}
\end{eqnarray}
\end{example}

\begin{remark}
The polar coordinates in \eqref{angemb1a}, \eqref{angemb1b} might be viewed as a geometric construction of the odd dimensional sphere $S^{2n+1} (\mu)$ with radius $\mu >0$ 
 of the following type: we consider the product of 
 the $n$--dimensional spheres  $S^n (2^{-1/2}\mu)$ in $n+1$--dimensional real subspaces 
 $\Im \zeta =0$ and $\Re\zeta =0$ plus a rotation with the angle $\theta$.
 \end{remark}

We propose an apparently new representation of the Kepler manifold by means of the co-dimension $2$ sub--manifold of $\C^{n+1}\setminus 0$ using  the angular variables  
\eqref{angemb1a}, \eqref{angemb1b}.

\begin{proposition}
The Kepler manifold $C$ is defined in the polar coordinates $r>0$, $(\varphi,\psi,\theta)$ in  
\eqref{angemb1a}, \eqref{angemb1b} by 
\begin{eqnarray}
 \{ (r,\varphi, \psi,\theta):\, \theta  = \frac{\pi}{4}, \ \ 
H(\varphi, \psi) =0\},  
\label{polar3}
\end{eqnarray}
where 
\begin{eqnarray}
H(\varphi, \psi) & := & e(\varphi)   \cdot s(\psi) \nonumber\\
& =  & 
\cos (\varphi_n -\psi_n) \prod_{j=1}^{n-1}(\sin \varphi_j \sin \psi_j) + 
\widetilde{H}(\varphi', \psi') 
\label{polar3a1}
\end{eqnarray}
with 
 $\varphi'= (\varphi_1, \ldots, \varphi_{n-1})$, $\psi'= (\psi_1, \ldots, \psi_{n-1})$ and 
\begin{eqnarray}
\widetilde{H}(\varphi', \psi') & = & \cos \varphi_1 \cos \psi_1 + \sum_{j=2}^{n-1} 
\cos \varphi_j \cos \psi_j  
 \prod_{\ell=1}^{j-1} \sin \varphi_\ell \sin\psi_\ell    .
\label{polar3b}
\end{eqnarray}
Clearly,   $\widetilde{C} $ is embedded in 
 $S^{2n+1}(\sqrt{2})$ by the equations  
\begin{eqnarray}
r = \sqrt{2}, \ \ \ \theta = \frac{\pi}{4},   \ \ \ 
H(\varphi, \psi) =0 
\label{polar3c}
\end{eqnarray}
and therefore 
\begin{eqnarray}
\rho & = & |x| = \frac{r}{\sqrt{2}}.
\label{polar3d}
\end{eqnarray}
Since $\widetilde{C} $ is compact 
we can find a finite covering $\ds \widetilde{C}= \bigcup_{k=1}^d\widetilde{U}_k$  of open charts $\widetilde{U}_k $, $k=1,\ldots, d$, which yields 
$$C= \bigcup_{k=1}^d ]0,+\infty[\times \widetilde{U}_k,$$ 
where each 
$\widetilde{U}_k $ is diffeomorphic to an open set  $U_k\subset  \R^{2n-1}$, with local coordinates $\Theta=\Theta^k=(\Theta^k_1, \ldots, \Theta^k_{2n-1})  \in U_k$. 
For every $k\in \{1, \ldots d\}$ the  two form 
$\partial\bar \partial  f  $, $f$ being a smooth function on $C$,  can be written in the (cylindric) coordinates 
$(\rho, \Theta) \in ]0, +\infty[\times U_k$  as follows 
\begin{eqnarray}
\partial\bar\partial f  &= & 
\sum_{\ell=1}^{2n-1} 
\rho  \left( \theta^\ell_0\partial_\rho^2 + \theta^\ell_1\rho^{-1}
 \partial_\rho L_1^{\ell} + \rho^{-1}b^\ell \partial_\rho + \rho^{-2} \widetilde{\Delta}_\ell \right)f (\rho, \Theta)  d\rho  \wedge  d\Theta_\ell
\nonumber\\
& + & \sum_{j,\ell=1}^{2n-1} 
\rho^2  \left( \theta^{j\ell}_0 \partial_\rho^2 + \theta^{j\ell}_1\rho^{-1}
 \partial_\rho L_1^{j\ell} + \rho^{-1}b^{j\ell} \partial_\rho + \rho^{-2} \widetilde{\Delta}^{j\ell} \right)f (\rho, \Theta)  d\Theta_j \wedge  d\Theta_\ell, 
\label{LB1}
\end{eqnarray}
where $\theta^\ell_0$, $\theta^{j\ell}_0 $, $\theta^\ell_1$, $\theta^{j\ell}_1$  
$b^\ell$, $b^{j\ell}$  are  real--valued real analytic functions on $U_j$,
 $L_1^\ell$, $L_1^{j\ell}$ are real 
 tangential vector fields to $\widetilde{C}$ with real analytic  coefficients 
while $\widetilde{\Delta}^{j\ell}$, 
$\widetilde{\Delta}^{j\ell}$ are second order linear analytic differential operators without zero order term  on $\widetilde{C}$ (i.e.,   $ \Delta^{j\ell} 1=
\widetilde{\Delta}^{j\ell} 1 =0$),   $j,\ell=1,\ldots,2n-1 $. 
In particular, if $f$ is constant on $\widetilde{C}$, i.e., $f= f(\rho)$, we have 
\begin{eqnarray}
\partial\bar\partial f  &= & 
\sum_{\ell=1}^{2n-1} 
\left( \theta^\ell_0 (\Theta) \rho f''(\rho) + b^\ell (\Theta) f' (\rho)\right)  d\rho  \wedge d\Theta_\ell
\nonumber\\
& + & \sum_{j,\ell=1}^{2n-1} 
\left(  \theta^{j\ell}_0(\Theta) \rho^2 f''(\rho) + b^{j\ell}(\Theta) \rho  
f'(\rho) \right)   d\Theta_j \wedge   d\Theta_\ell. 
\label{LB1a}
\end{eqnarray}
\end{proposition}

 {\it Proof}: 
 We start by recalling the representation of 
 $ \partial\bar\partial f$ in $\C^{n+1}$ by identifying $\C^{n+1}= \R^{2n+2}$ by the canonical 
 complex structure $J_0$ for $n\geq 1$. Recall that if $n=0$ we have  
 $$ \partial\bar\partial f = 2i (f_{xx} + f_{yy}) dx \wedge dy,$$ 
 with $z= x+yi \in \C$.

 \begin{lemma} Let $n\geq 1$ and $f\in C^\infty (\C^{n+1}\setminus 0) = 
 C^\infty (\R^{2n+2}\setminus 0)$ with coordinates $\zeta = \xi + i\eta$. 
 Then we can write $\partial \bar\partial f$ in the standard polar coordinates 
 in $\R^{2n+2}\setminus 0$:  
 \begin{eqnarray}
\partial\bar\partial f  &= & \frac{i}{2}
\sum_{\ell=1}^{n+1} \partial_{\zeta_\ell}\partial_{\bar \zeta_\ell} f  d\zeta_\ell \wedge d\bar\zeta_\ell 
\nonumber\\
= \sum_{\ell=1}^{n+1} ( f_{\xi_\ell\xi_\ell} + f_{\eta_\ell \eta_\ell})  d\xi_\ell \wedge d\eta_\ell 
\label{can2f1}
\end{eqnarray}
Moreover, in the the standard polar coordinates 
 in $\R^{2n+2}\setminus 0$ defined by \eqref{polar1can} we can write 
 \begin{eqnarray}
\partial\bar\partial f  &= & 
\rho \sum_{\ell=1}^{2n+1}    \Omega^\ell[f]  
  d\rho  \wedge  d\Theta_\ell + \rho^2 \sum_{j,\ell=1}^{2n+1} 
\Omega^{j\ell}[f] d\Theta_j \wedge  d\Theta_\ell, 
\label{can2f2}
\end{eqnarray}
with 
\begin{eqnarray}
\Omega^\ell [f]  &  =  &   
 E^\ell_0(\Theta) \partial_\rho^2 f + E^\ell_1(\Theta) \rho^{-1}
 \partial_\rho L_1^{\ell;E}(\Theta, \partial_\Theta)f \nonumber\\
 & + & \rho^{-1}e^\ell(\Theta) \partial_\rho f + \rho^{-2} \widetilde{\Delta}_\ell^E(\Theta, \partial_\Theta) f 
 \label{can2f2a}\\
\Omega^{j\ell}[f]  & = &   
E^{j\ell}_0(\Theta) \partial_\rho^2 f + E^{j\ell}_1(\Theta)\rho^{-1}
 \partial_\rho L_1^{j\ell;E}(\Theta,\partial_\Theta) f 
 \nonumber\\
 & + &  \rho^{-1}e^{j\ell}(\Theta) \partial_\rho  f+ \rho^{-2} \widetilde{\Delta}^{j\ell;E}(\Theta,\partial_\Theta)  f 
\label{can2f2b}
\end{eqnarray}
where $E^\ell_0$, $E^{j\ell}_0 $, $E^\ell_1$, $E^{j\ell}_1$  
$e^\ell$, $e^{j\ell}$  are  real--valued real analytic functions of $\theta$,
 $L_1^{\ell;E}$, $L_1^{j\ell;E}$ are real 
 tangential vector fields to $S^{2n+1}$ with real analytic  coefficients 
while $\widetilde{\Delta}^{j\ell;E}$, 
$\widetilde{\Delta}^{j\ell;E}$ are second order linear analytic differential operators without zero order term  on $S^{2n+1}$ (i.e.,   $ \Delta^{j\ell;E} (1)=
\widetilde{\Delta}^{j\ell;E} (1) =0$),   $j,\ell=1,\ldots,2n+1 $.  
 \end{lemma}

 \textit{Proof}. The first identity \eqref{can2f1} is immediate. Next, by the 
 standard calculus on two forms and change of the variables, we obtain
 \begin{eqnarray}
\partial\bar\partial f  &= & \sum_{\ell=1}^{n+1} ( f_{\xi_\ell\xi_\ell} + f_{\eta_\ell \eta_\ell})  d\xi_\ell \wedge d\eta_\ell 
\nonumber\\
& = & \sum_{j=1}^{n+1} ( f_{\xi_j\xi_j} + f_{\eta_j \eta_j}) \sum_{\ell=1}^{2n+1} 
(\partial_\rho \xi_j \partial_{\Theta_\ell}\eta_j - \partial_\rho \eta_j \partial_{\Theta_\ell}\xi_j) d\rho \wedge d\Theta_\ell \nonumber\\
& + &  \sum_{j=1}^{n+1} ( f_{\xi_j\xi_j} + f_{\eta_j \eta_j}) \sum_{k,\ell=1}^{2n+1} 
(\partial_{\Theta_k} \xi_j \partial_{\Theta_\ell}\eta_j - \partial_{\Theta_k}\eta_j \partial_{\Theta_\ell}\xi_j) d\Theta_k \wedge d\Theta_\ell\nonumber\\
& = & \sum_{\ell=1}^{2n+1} \left( \sum_{j=1}^{n+1} ( f_{\xi_j\xi_j} + f_{\eta_j \eta_j})
(\partial_\rho \xi_j \partial_{\Theta_\ell}\eta_j - \partial_\rho \eta_j \partial_{\Theta_\ell}\xi_j)
\right) d\rho \wedge d\Theta_{\ell}\nonumber\\
& + & 
\sum_{k,\ell=1}^{2n+1} \left( \sum_{j=1}^{n+1} ( f_{\xi_j\xi_j} + f_{\eta_j \eta_j})
(\partial_{\Theta_k} \xi_j \partial_{\Theta_\ell}\eta_j - \partial_{\Theta_k} \eta_j \partial_{\Theta_\ell}\xi_j)
\right) d\Theta_k \wedge d\Theta_{\ell}.
\label{can2f3}
\end{eqnarray}


We conclude by writing 
$f_{\xi_j\xi_j} + f_{\eta_j \eta_j}$ in the polar coordinates $(\rho, \Theta)$. 
\vspace{10pt}

Set $\Phi = (\varphi, \psi)\in I\times J$, where $I$ (respectively $J$) is 
as in  \eqref{defpol1} (respectively, \eqref{defpol2}).
We can apply the implicit function theorem in the global coordinates on $\widetilde{C}$ 
is applicable outside  the  singular set $\widetilde{C}_S(I\times J)$ (in the fixed angular variables) defined by the system 
\begin{eqnarray}
H(\Phi) & = & 0, \qquad d_\Phi H(\Phi) = 0.
\label{polar3a}
\end{eqnarray}
In fact, choose 
$V\subset\subset I\times J\setminus  \widetilde{C}_S(I\times J)$, then 
by the implicit function theorem, there exists $s\in \{1,\ldots, n\}$ such that $H(\Phi)=0$ 
on $V$ defines 
$\varphi_s = \varphi_s( \varphi^s, \psi)$ or $\psi_s = \psi_s (\varphi, \psi^s)$, 
where 
$$\varphi^s = (\varphi_1, \ldots,\varphi_{s-1},\varphi_{s+1}, \ldots, \varphi_n), \qquad  
\varphi^s = (\varphi_1, \ldots,\varphi_{s-1},\varphi_{s+1}, \ldots, \varphi_n).
$$ 
with $\Theta =(\varphi^s, \psi)$ or $\Theta =(\varphi, \psi^s)$ belonging to some open
 $V' \subset \R^{2n-1}$. 
Clearly $V'$ defines a chart and in view of the compactness of $S^{2n+1}(r)$ we can choose a finite number by varying $I$ and $J$ in the definitions \eqref{defpol1}, \eqref{defpol2}.

It is well known by the calculus on manifolds 
that 
\begin{eqnarray}
\partial_{z_j} \partial_{\overline{z_\ell} } f & = &  
\sum_{r=1}^{2n-1} 
( \Omega^{j\ell}_r(\Theta) \partial_\rho^2 + \rho^{-1}
 \partial_\rho T^{j\ell}_r + \rho^{-1} \tilde{b}^{j\ell}_r \partial_\rho + 
 \rho^{-2}  D^{j\ell}_r ) f   
\nonumber \\
 & + & \sum_{r,s=1}^{2n-1} 
(\Omega^{j\ell}_{rs} \partial_\rho^2 + \rho^{-1}
 \partial_\rho T^{j\ell}_{rs} + \rho^{-1}\tilde{b}^{j\ell}_{rs} \partial_\rho + 
 \rho^{-2} D^{j\ell}_{rs} ) f   
\label{GenLB1}
\end{eqnarray}  
and 
\begin{eqnarray}
dz_j\wedge d\bar z_{\ell}    
& = & \rho \sum_{p=1}^{2n-1} \Gamma^{j\ell}_p(\Theta) d\rho \wedge d\Theta_p \nonumber\\
& + & \rho^2 \sum_{p,q=1}^{2n-1}
\Gamma^{j\ell}_{pq}(\Theta) d\Theta_p \wedge d\Theta_q 
\label{GenLB1a}
\end{eqnarray}
 where  
 \begin{eqnarray}
 T^{j\ell}_r & = & \sum_{q=1}^{2n-1} T^{j\ell}_{r;q}(\Theta) \partial_{\Theta_q}
 \label{T1}\\
 T^{j\ell}_{rs}& = & \sum_{q=1}^{2n-1} T^{j\ell}_{rs;q}(\Theta) \partial_{\Theta_q}
 \label{T2}
 \end{eqnarray}
 are  tangential vector field to $\widetilde{C}\bigcap \tilde{U}_j$ with real analytic  coefficients 
while
\begin{eqnarray}
 D^{j\ell}_r  & = & \sum_{p,q=1}^{2n-1}  D^{j\ell}_{r;pq}(\Theta) \partial_{\Theta_p}\partial_{\Theta_q} \nonumber\\
 &  + &  \sum_{q=1}^{2n-1} 
 H^{j\ell}_{r;q}(\Theta) ) \partial_{\Theta_q}
 \label{D1}\\
 D^{j\ell}_{rs}& = & 
 \sum_{p,q=1}^{2n-1} D^{j\ell}_{rs;pq}(\Theta) \partial_{\Theta_p}\partial_{\Theta_q}
 \nonumber\\
 & + & \sum_{q=1}^{2n-1} 
  H^{j\ell}_{rs;q}(\Theta) \partial_{\Theta_q} 
 \label{D2}
 \end{eqnarray}
are second order linear analytic differential operator on  $\widetilde{C}\bigcap \tilde{U}_j$  without zero order terms.   
Clearly, \eqref{obstr1} and  \eqref{GenLB1} yield  \eqref{LB1}. The proof is complete. 
\vspace{10pt}

\vspace{10pt}

We observe  that the definition of 
$T_m(z)$ 
implies that 
\begin{eqnarray}
\log T_m(z) &= & n \log (2m) +  \log \left(  \frac{T_m(z)}{2m^n} \right)\nonumber\\
& = & n \log (2m) + \log F( m|x|)\nonumber\\
& = &  n \log (2m) +
\log \left( \sum_{j=0}^{n-2}  \frac{a_j}{(m|x|)^j}   + \Phi(m|x|) + \Psi(m|x|)\right).
\label{obstr1}
\end{eqnarray}
Set 
\begin{eqnarray}
F(\rho) &= & \sum_{j=0}^{n-2}  \frac{a_j}{\rho^j}   + \Phi(\rho) + \Psi(\rho).
\label{obstr1a}
\end{eqnarray}
Clearly we can rewrite \eqref{obstr1} on $C= ]0,+\infty[\times \tilde{C}$ as follows:  
\begin{eqnarray}
\log T_m(z) 
& = &  n \log (2m) +
\log F(m |x|) =n \log (2m) +
\log F(m \rho).  
\label{obstr1b}
\end{eqnarray}
We note that \eqref{obstr1} implies for      every fixed $\rho >0$ the function 
 $\log T_m(z)$ is constant $C= \{ \rho \}\times \tilde{C}$. 
 The representaion formula   \eqref{obstr1b} yileds 
\begin{eqnarray} 
 {\mathcal E}_m|_C & = & \frac{i}{2m} \partial \bar \partial ( \log (F (m |x|) |_C.
\label{obstr1c}
\end{eqnarray}

We show the main result of the present section.
\vspace{10pt}

\begin{theorem}   
There exist   $2(n+1)^2$  real analytic functions $\sigma_{k\ell}  (\Theta)$, $
\tau_{k\ell} (\Theta)$, 
$\overline{\sigma_{k\ell}} = \sigma_{\ell k}$, 
$\overline{\tau_{k\ell}} = \sigma_{\ell k}$, $k,\ell =1,\ldots, n+1$,   
 defined on $\widetilde{C}$ such the $2$--form 
the obstruction term  
\begin{eqnarray} 
 {\mathcal E}_m|_{C_j} (\rho, \Phi) 
 & := & {\mathcal E}_m|_{]0,+\infty[\times (\tilde{C} \bigcap U_j} (\rho, \Phi) \nonumber\\
 & = & \sum_{q=1}^{2n-1}  {\mathcal E}_m^{q}(\rho,\Theta)
  d\rho  \wedge d\Theta_q
\nonumber\\
& + & \sum_{p,q=1}^{2n-1} 
{\mathcal E}_m^{pq}(\rho,\Theta)
  d\Theta_p \wedge d\Theta_q 
\label{eras2}
\end{eqnarray}
where
 \begin{eqnarray}
{\mathcal E}_m^{q}(\rho,\Theta) & = & 
\theta^q_0 (\Theta) m \rho \frac{F''(m\rho) F(m\rho) -(F'(m\rho))^2}{F^2(m\rho)} 
 + b_q(\Theta) \frac{F'(m\rho)}{F(m\rho)} \nonumber\\
 & = &  \frac{1}{m^2\rho^2} \left( 
 \theta^q_0 (\Theta) 
 \sum_{s=0}^N \frac{p_s}{\rho^s m^s} 
 + b_q (\Theta) 
 \sum_{s=0}^N \frac{q_s}{\rho^s m^s} \right) 
  \nonumber\\ 
  & + & \theta^q_0 (\Theta) R^{1}_N(m\rho) + b_q(\Theta) R^{2}_N (m\rho) 
\label{eras2a }\\
{\mathcal E}_m^{pq}(\rho,\Theta)
& = &  
 \theta_0^{pq}(\Theta) m\rho^2 \frac{F''(m\rho) F(m\rho) -(F'(m\rho))^2}{F^2(m\rho)}
  +  b_{pq}(\Theta) \rho  
\frac{F'(m\rho)}{F(m\rho)}  \nonumber\\
& = & \frac{1}{m^2\rho} \left( \theta^{pq}_0 (\Theta) \sum_{s=0}^N \frac{p_s}{\rho^s m^s} 
 + b_{pq} (\Theta) \sum_{s=0}^N \frac{p_s}{\rho^s m^s} \right) 
 \nonumber\\
  & + & 
 \theta^{pq}_0 (\Theta) \rho R^{1}_N(m\rho) + b_{pq}(\Theta) \rho R^{2}_N (m\rho) 
\label{eras2b}
\end{eqnarray}
for all $N\in \N$, where the real constants $p_s$, $q_s$ depend on   
$a_1, \ldots, a_{n-2}$ by explicit formulas  while  
$\theta^{q}_0 (\Theta)$, $b_{q}(\Theta)$,  
$\theta^{pq}_0 (\Theta)$, $b_{pq}(\Theta)$ are real-valued real analytic functions 
on $\tilde{C}$.  
The reminders  $R^{\mu}_{N}$, $\mu=1,2$  
satisfy  the following estimates: for every 
$\delta_0 >0$ 
one can find positive constants $A_{\mu}$, $B_{\mu}$, $\mu=1,2$, 
such that  
\begin{eqnarray}
|\partial_x ^\alpha (R_{N}^\mu ( m|x|))|  & \leq  & A_{\mu}^{N+1} B_{\mu}^{|\alpha|} N! \alpha! m^{-N-1}|x|^{-N-|\alpha|-3}, 
\label{eras3}
\end{eqnarray}
for all $ N\in \Z_+$, $\alpha \in \Z_+^n$, $|x|\geq \delta_0$, 
$m\geq 1$. In particular, there exist positive  constants $c_0, \mu_0, A_0$ 
such that if we choose $\ds N = N(m) = 1 + [e^{c_0 m}]$ 
then $\widetilde{R^\mu}_m (m\rho ) := R^\mu_{N(m)} (m\rho)$ are 
exponentially small in the following sense: 
\begin{eqnarray}
|\partial_x ^\alpha (\widetilde{R^\mu}_m ( m|x|))|  & \leq & A_0B_\mu^{|\alpha|}\alpha! e^{- c_0 m |x| }, 
\label{eras4}
\end{eqnarray}
for all $\alpha \in \Z_+^n$, $|x|\geq \delta_0$, $m\geq 1$. 
Finally, we can summarize the asymptotic estimates above as follows:
 every $N\in \N$, $\delta >0$ one can find $C>0$ such that  
\begin{eqnarray} 
{\mathcal E}_m^{q} (\rho, \Theta) & = & 
\frac{1}{\rho^2 m^2} \sum_{s=1}^{N} \frac{\varkappa_s^{q}(\Theta) }{m^s\rho^{s}}  + R_N^q(m\rho, \Theta)
\label{obstrest1}\\
{\mathcal E}_m^{pq} (\rho, \Theta) & = & 
\frac{1}{\rho m^2} \sum_{s=1}^{N} \frac{\varkappa_s^{pq}(\Theta) }{m^s\rho^{s}}  + \rho R_N^{pq}(m\rho, \Theta), 
\label{obstrest1a}
\end{eqnarray}
where 
\begin{eqnarray}
|\partial_\rho^\beta \partial^\alpha_\Theta (R_N^{q} ( \rho, \Theta;m))|  & \leq  & C^{N+1 +|\alpha|} N! \alpha! m^{-N-3}\rho^{-N-3-\beta}, 
\label{obstrest2}\\
|\partial_\rho^\beta \partial^\alpha_\Theta (R_N^{pq} ( \rho, \Theta;m))|  & \leq  & C^{N+1 +|\alpha|} N! \alpha! m^{-N-3}\rho^{-N-2-\beta}, 
\label{obstrest2a}
\end{eqnarray}
for $m\geq 1$, $|x|\geq \delta$, $\beta \in \Z_+$, $\alpha \in \Z_+^{2n-1}$, 
$k,\ell =1,\ldots, n+1$, $p,q=1,\ldots 2n-1$. 
\end{theorem}
 
\textit{Proof}. 
Fix $m\in N$. The first step of the proof consists in applying \eqref{LB1a} for 
$$ f(\rho) = \frac{i}{2m}\log (F(m\rho)).$$ 
We note that by \eqref{LB1a} for $f(\rho) = \frac{i}{2m}\log (F(m\rho))$ and the identities  
\begin{eqnarray}
\partial_\rho (\log (F(m\rho) ) & = & m\frac{F'(m\rho) }{F(m\rho)},   
\label{pobst1}\\
 \partial_\rho^2 (\log (F(m\rho) ) & = & m^2\frac{F''(m\rho)F(m\rho) -(F'(m\rho))^2}{F^2(m\rho)}, 
\label{pobstr2}
\end{eqnarray}
we get 
\begin{eqnarray}
{\mathcal E}_m^{q}(\rho,\Theta) & = & \frac{i}{2}
\frac{ \Gamma^{q} ( \Theta, m\rho) }{ F^2 (m\rho) } \label{pobstr3}\\
{\mathcal E}_m^{pq}(\rho,\Theta) & = & \frac{im\rho}{2} \frac{ \Gamma^{pq} (\Theta, m\rho) }{ F^2 (m\rho) } \label{pobstr3a}
\end{eqnarray}
where 
\begin{eqnarray}
\Gamma^q(\Theta, y) & = &  \theta_0^{q}(\Theta) y( F''(y)F(y) -(F'(y))^2) + 
b^q(\Theta) F'(y) F(y), 
\label{obstr2c}\\
\Gamma^{pq}(\Theta, y) & = &  \theta_0^{pq}(\Theta) y( F''(y)F(y) -(F'(y))^2) + 
b^{pq}(\Theta) F'(y) F(y). 
\label{obstr2d}
\end{eqnarray}

By the  asymptotic expansion for $F$ we have 
\begin{eqnarray}
F'(y) &= & -\sum_{j=1}^{n-2}  \frac{j b_j}{y^{j+1}}  + \Phi'(y) + \Psi'(y)
\label{rep1D1}\\
F''(y) &= & \sum_{j=1}^{n-2}  \frac{j(j+1) b_j}{y^{j+2}}  + \Phi''(y) + \Psi''(y)
\label{rep1D2}
\end{eqnarray}

Straightforward calculations of $F''(y) F(y)$, $(F'(y))^2$, $F'(y)F(y)/y$ and \eqref{obstr2c}  
lead to 
\begin{eqnarray}
\Gamma^q(\Theta, y) & = &  \theta_0^{q}(\Theta) y( F''(y)F(y) -(F'(y))^2) + 
b^q(\Theta) F'(y) F(y), 
\label{obstr2c1}\\
\Gamma^{pq}(\Theta, y) & = &  \theta_0^{pq}(\Theta) y( F''(y)F(y) -(F'(y))^2) + 
b^{pq}(\Theta) F'(y) F(y). 
\label{obstr2d1}
\end{eqnarray}

\begin{eqnarray}
\Gamma^q(\Theta, y) & = &  
\sum_{j=1}^{2n-4} \frac{\mu_{j;0}^q(\Theta)}{y^{j+2} } +   
 \sum_{j=1}^{2n-4} \frac{\nu^q_{j;0}(\Theta) }{y^{j+2}}\nonumber \\
& + &  \theta_0^q (\Theta) E_1(y) + b^q(\Theta) E_2(y)  
 \label{rep2a}\\
 \Gamma^{pq}(\Theta, y) & = &  
\sum_{j=1}^{2n-4} \frac{\mu_{j;0}^{pq}(\Theta)}{y^{j+2} } +   
 \sum_{j=1}^{2n-4} \frac{\nu^{pq}_{j;0} (\Theta)  }{y^{j+2}}\nonumber \\
& + &  \theta_0^{pq} (\Theta) E_1(y) + b^{pq}(\Theta)(\Theta) E_2(y)  
\label{rep2b}
\end{eqnarray}
where 
\begin{eqnarray}
\mu_{1;0}^q(\Theta) & = & \theta_0^{q}(\Theta)  b_1 b_0 \label{crep1a}\\
\mu_{j;0}^q(\Theta)  &= &  \theta_0^{q}(\Theta) \left(
b_0 b_j + \sum_{\ell=1}^{j-1} \ell (\ell+1) b_\ell b_{j-\ell}   
- \sum_{\ell=1}^{j-1} \ell (j -\ell) b_\ell  b_{j-\ell}\right) \nonumber\\ 
& = & \theta_0^{q}(\Theta) \left(
j(j+1) b_0 b_j + \sum_{\ell=1}^{j-1} \ell (2\ell+1-j) b_\ell b_{j-\ell}\right)  
\label{crep2a} \\
\nu_{1;0}^q(\Theta) & = & -b^{q}(\Theta) b_1 b_0 \label{crep1b}\\
\nu_{j;0}^q(\Theta) & = & -\theta_0^{q}(\Theta) 
(jb_0 b_j + \sum_{\ell=1}^{j-1} \ell  b_\ell b_{j-\ell})  
\label{crep2b}\\
\mu_{1;0}^{pq}(\Theta) & = & \theta_0^{q}(\Theta)  b_1 b_0 \label{crep1aa}\\
\mu_{j;0}^{pq}(\Theta)  &= &  \theta_0^{pq}(\Theta) \left(
b_0 b_j + \sum_{\ell=1}^{j-1} \ell (\ell+1) b_\ell b_{j-\ell}   
- \sum_{\ell=1}^{j-1} \ell (j -\ell) b_\ell  b_{j-\ell}\right) \nonumber\\ 
& = & \theta_0^{pq}(\Theta) \left(
j(j+1) b_0 b_j + \sum_{\ell=1}^{j-1} \ell (2\ell+1-j) b_\ell b_{j-\ell}\right)  
\label{crep2aa} \\
\nu_{1;0}^{pq}(\Theta) & = & -b^{pq}(\Theta) b_1 b_0 \label{crep1bb}\\
\nu_{j;0}^{pq}(\Theta) & = & -\theta_0^{pq}(\Theta) 
(jb_0 b_j + \sum_{\ell=1}^{j-1} \ell  b_\ell b_{j-\ell})  
\label{crep2bb}
\end{eqnarray}
for $ j=2,\ldots, 2n-4$ 
and 
\begin{eqnarray}
E_1(y) & = & (\Phi (y) +\Psi (y)) \sum_{j=1}^{n-2}  \frac{j(j+1) b_j}{y^{j+2}} + 
(\Phi'' (y) +\Psi'' (y)) \sum_{j=0}^{n-2}  \frac{ b_j}{y^{j}} 
 \nonumber\\ 
& - & 2 (\Phi' (y) +\Psi' (y)) \sum_{j=1}^{n-2}  \frac{j b_j}{y^{j+1}}\nonumber\\
& + & (\Phi (y) +\Psi (y)) (\Phi'' (y) +\Psi'' (y))
\label{erep1a}\\
E_2(y) & = & -(\Phi (y) +\Psi (y)) \sum_{j=1}^{n-2}  \frac{j b_j}{y^{j+2}} - 
(\Phi' (y) +\Psi' (y)) \sum_{j=0}^{n-2}  \frac{ b_j}{y^{j+1}}\nonumber\\
& - & (\Phi (y) +\Psi (y))(\Phi' (y) +\Psi' (y)) 
\label{erep1b}
\end{eqnarray}

Next, we show an auxiliary  assertion.

\begin{lemma}
There exists $K>0$ such that the function $G(y):=F^{-2}(y)$ is uniformly analytic function for $y\geq 2$ satisfying 
\begin{eqnarray}
\frac{1}{F^2(y)} & = & \frac{1}{b_0^2} + \sum_{j=1}^\infty \frac{\beta_j}{y^j}, 
\label{den1}
\end{eqnarray}
where $\beta_j \in \R$, $j=1,2, \ldots $ are defined recursively  
\begin{eqnarray}
\beta_j  & = & -2\frac{b_j}{b_0} + 
\widetilde{\beta_j}(b_0, \ldots, b_{j-1}), \qquad
j=1,2, \ldots 
\label{den2rec}
\end{eqnarray}
 and  satisfy 
\begin{eqnarray}
\limsup_{j\geq 1}\sqrt[j]{|\beta_j|} & < & \frac{1}{K}.
\label{den2}
\end{eqnarray}
Moreover, 
\begin{eqnarray}
\frac{1}{F^2(m|x|)} & = & \frac{1}{b_0^2} + \sum_{j=1}^\infty \frac{\beta_j}{(m|x|)^j}, 
\label{den3}
\end{eqnarray}
 is uniformly analytic  for $|x| \geq K$, uniformly with respect to $m \geq 1$ and the remainder    
\begin{eqnarray}
\widetilde{E}_N (y)& = & \sum_{j=N}^\infty \frac{\beta_j}{y^j}
\label{den4}
\end{eqnarray}
satisfies, for some positive numbers 
$\widetilde{A}_0$,  $\widetilde{B}_0$,  $\widetilde{A}$, 
  $\widetilde{B}$, 
 the following analytic--Gevrey combinatorial  
estimates 
\begin{eqnarray}
|(\frac{d}{dy}) ^k (\widetilde{E_N}(y ))|  & \leq & \widetilde{A}_0^{N+1} 
\widetilde{B}_0^k N! k! y^{-N-k-1}, 
\label{den5a}
\end{eqnarray} 
for $y \geq K$,  $k, N\in \Z_+$,  and 
\begin{eqnarray}
|\partial_x ^\alpha (\widetilde{E_N}( m|x|))|  & \leq & \widetilde{A}^{N+1} \widetilde{B}^{\alpha} N! \alpha! m^{-N-1}|x|^{-N-\alpha-1}, 
\label{den5}
\end{eqnarray} 
\end{lemma} 
for all  $\alpha \in \Z_+^n $, $N\in \Z_+$,  $|x|\geq K$, $m\geq 1$. 

\textit{Proof}. We recall that 
\begin{eqnarray}
\frac{1}{(1+\tau)^2} & = & -\frac{d}{d\tau}
\left( \frac{1}{1+\tau }\right)  =   \sum_{j=1}^\infty (-1)^j j\tau^{j-1} 
\label{den6}
\end{eqnarray}
provided $|\tau| <1$. 
Set 
\begin{eqnarray}
\varkappa(y) & = & \frac{1}{b_0} \sum_{j=1}^{n-2}  \frac{b_j}{y^{j}}  + \frac{1}{b_0}\Phi(y) + \frac{1}{b_0}\Psi(y).
\label{den7}
\end{eqnarray}
Clearly for every $\veps \in ]0,1[$ we can find $K=K_\veps >0$ such that 
\begin{eqnarray}
\sup_{y\geq K} |\varkappa(y)| & < &\veps.
\label{den8}
\end{eqnarray}
Therefore, by \eqref{den7} and \eqref{den8} we readily obtain 
the following representation by means of  convergent Neumann series  
\begin{eqnarray}
\frac{1}{F^2(y)} & = & \frac{1}{b_0^2}
\frac{1}{(1+\varkappa (y))^2} \nonumber\\
&  = &  \frac{1}{b_0^2}  + \sum_{j=1}^\infty \frac{(-1)^{j+1} (j+1) }{b_0^2} 
(\varkappa (y))^j. 
\label{den9}
\end{eqnarray}
Hence, 
\begin{eqnarray}
\beta_j & = & \frac{1}{j!} 
(\frac{d}{d y})^j (\sum_{s=1}^\infty \frac{(-1)^{s+1} (s+1) }{b_0^2} 
(\varkappa (y))^j)| _{y=0}.   
\label{den9a}
\end{eqnarray}

We complete the proof of the lemma by plugging \eqref{den7} in \eqref{den9} and 
applying uniform analytic estimates as in  \cite{cgrjfa},  \cite{cgrbir}.
\vspace{10pt}

Next, by the representation formulas derived above, we obtain that  
\begin{eqnarray}
\frac{ \Gamma^q (\Theta, y )}{ F^2(y)} & = & 
\sum_{j=1}^{N} \frac{H_j^q(\Theta,y)}{y^{j+2} } +   
\sum_{j=1}^{N} \frac{G_j^q(\Theta,y }{y^{j+2} } 
\nonumber \\
& + & \theta_0^q(\Theta)
\widetilde{R^1_N}(\Theta, y) + 
b^q(\Theta)\widetilde{R^2_N}(y)
\label{den10a}\\
\frac{ \Gamma^{pq} (\Theta, y )}{ F^2(y)} & = & 
 \sum_{j=1}^{N} \frac{H_j^{pq}(\Theta,y)}{y^{j+2} } +   
 \sum_{j=1}^{N} \frac{G_j^{pq}(\Theta,y }{y^{j+2} } 
\nonumber \\
& + & \theta_0^{pq}(\Theta)
\widetilde{R^1_N}(\Theta, y) + 
b^{pq}(\Theta)\widetilde{R^2_N}(y)
\label{den10b}
\end{eqnarray}
where
\begin{eqnarray}
 H_j^q(\Theta,y) & = & \theta^q_0(\Theta) \left(
\frac{ \mu_j }{b_0^2}  +  
\sum_{\ell = 1 }^{ \min \{2n-4, j-2 \} }  
\mu_\ell \beta_{j-\ell} \right)
\label{den11a}\\
 H_j^{pq}(\Theta,y) & = & \theta^{pq}_0(\Theta) \left(
\frac{ \mu_j }{b_0^2}  +  
\sum_{\ell = 1 }^{ \min \{2n-4, j-2 \} }  
\mu_\ell \beta_{j-\ell} \right)
\label{den11b}\\
G^q_j(\Theta,y) 
 & = &  \theta^q_0(\Theta) \left(
\frac{ \nu_j }{b_0^2}  +  
\sum_{\ell = 1}^{ 
 \min  \{ 2n-4, j-1 \} } 
\nu_\ell \beta_{j-\ell} \right)
\label{den12a}\\
G^{pq}_j(\Theta,y) 
 & = &  \theta^{pq}_0(\Theta) \left(
\frac{ \nu_j }{b_0^2}  +  
\sum_{\ell = 1}^{ 
 \min  \{ 2n-4, j-1 \} } 
\nu_\ell \beta_{j-\ell} \right)
\label{den12b}
\end{eqnarray}
 and 
\begin{eqnarray}
\widetilde{R^1_N}(y)
 & = &   
\sum_{\ell = 1}^{ 2n-4 } 
\sum_{N+1 -\ell  \leq j   \leq N}  
\frac{ \mu_\ell \beta_{j} }{y^{\ell +j +2 } } 
\nonumber \\
& + & 
\frac{E_1(y)}{F^2(y)}              +  
\widetilde{E_N}(y ) \sum_{j=1}^{2n-4} \frac{\mu_j}{y^{j+2} }, 
\label{den13}
\end{eqnarray}
\begin{eqnarray}
\widetilde{R^2_N}(y)
 & = &   
\sum_{\ell = 1}^{ 2n-4 } 
\sum_{N+1 -\ell  \leq j   \leq N}  
\frac{ \nu_\ell \beta_{j} }{y^{\ell +j +2 } } 
\nonumber \\
& + & 
\frac{E_2(y)}{F^2(y)}  +  
\widetilde{E_N}(y ) \sum_{j=1}^{2n-4} \frac{\nu_j}{y^{j+2} }
\label{den14}
\end{eqnarray}

We conclude the proof of the estimates by straightforward applications of 
the functional--analytic arguments in \cite{cgrjfa}, \cite{cgrbir} for showing simultaneosly uniform holomorphic  extensions  and exponential 
decay on infinity. 

\begin{remark}
We point out that  \eqref{eras4} implies  that for a given $m\gg 1$, the optimla choice of $N= 
N(m)$ for the  truncated asymptotic expansion  is given by 
\begin{equation}
\sigma(\Phi)  \sum_{j=1}^{N(m)}  \frac{p_j }{|x|^{j+2}} m^{-j} + 
\tau(\Phi)  \sum_{j=1}^{N(m)}  \frac{q_j }{|x|^{j+2}} m^{-j}.  
\label{eras5}
\end{equation}  
since the remainder is exponentially small for $m  
\to  +\infty $. 
We note that  similar uniform exponential decay estimates are shown 
in the framework of analytic--Gevrey pseudodifferential operators, e.g.,  cf. 
\cite{Sj1}, \cite{cgrjfa}, \cite{cgrbir} and the references therein. 
\end{remark}

\vskip 0.3cm

\noindent
{\bf Acknowledgements}. 
The authors thank  Cornelis Van Der Mee for his precious  comments and discussions on various 
aspects  in Quantum Mechancis.   
Thanks are also due  to Ivailo Mladenov and Vassil Tsanov 
for useful discussions and providing references related to the subject of the paper.

\small{}

\end{document}